\documentclass[3p]{elsarticle}
\usepackage{amssymb}
 \usepackage{amsthm}

\usepackage{amssymb,amsmath,amsthm,amsxtra,amsfonts}
\newcommand{\bbC}{{\mathbb{C}}}
\newcommand{\bbD}{{\mathbb{D}}}

\newcommand{\bbN}{{\mathbb{N}}}

\newcommand{\bbR}{{\mathbb{R}}}

\newcommand{\bdone}{{\boldsymbol{1}}}
\newcommand{\bdnot}{{\boldsymbol{0}}}

\newcommand{\llangle}{\left\langle\!\!\left\langle}
\newcommand{\rrangle}{\right\rangle\!\!\right\rangle}
\newcommand{\lla}{\left\langle\!\left\langle}
\newcommand{\rra}{\right\rangle\!\right\rangle}
\newcommand{\dt}{\frac{d\theta}{2\pi}}

\newcommand{\tr}{\text{\rm{Tr}}}

\newcommand{\beq}{\begin{equation}}
\newcommand{\eeq}{\end{equation}}
\newcommand{\ba}{\begin{align*}}
\newcommand{\ea}{\end{align*}}

\newcommand{\norm}[1]{\lVert#1\rVert}

\DeclareMathOperator{\real}{Re}
\DeclareMathOperator{\imag}{Im}

\DeclareMathOperator*{\esssup}{ess\,supp}

\DeclareMathOperator*{\res}{res}

\allowdisplaybreaks
\numberwithin{equation}{section}
\newtheorem{theorem}{Theorem}
\newtheorem{lemma}{Lemma}
\newtheorem{corollary}{Corollary}

\theoremstyle{remark}
\newtheorem*{remarks}{Remarks}
\begin{document}

\begin{frontmatter}



\title{Szeg\H{o} asymptotics for matrix-valued measures \\ with countably many bound states}


\author{Rostyslav Kozhan}
\ead{rostysla@caltech.edu}
\address{California Institute of Technology\\
Department of Mathematics 253-37\\
Pasadena, CA 91125, USA}

\begin{abstract}
Let $\mu$ be a matrix-valued measure with the essential spectrum a single interval and countably many point masses outside of it. Under the assumption that the absolutely continuous part of $\mu$ satisfies Szeg\H{o}'s condition and the point masses satisfy a Blaschke-type condition, we obtain the asymptotic behavior of the orthonormal polynomials on and off the support of the measure.

%Let $\mu$ be a matrix-valued measure with $\esssup\mu=[-2,2]$ and countably many point masses $\{E_j\}_{j=1}^\infty$ outside $[-2,2]$. Assume the Szeg\H{o} condition $\int_{-2}^2 (4-x^2)^{-1/2}\log\det \mu'_{a.c.}(x)dx>-\infty$ holds, and $\{E_j\}_{j=1}^\infty$ satisfy the Blaschke-type condition $\sum_{j=1}^\infty (|E_j|-2)^{1/2}<\infty$. Then Szeg\H{o} asymptotics for the orthogonal polynomials holds.

\noindent The result generalizes the scalar analogue of Peherstorfer--Yuditskii \cite{PY} %and Killip--Simon \cite{KS},
and the matrix-valued result of Aptekarev--Nikishin \cite{AN}, which handles only a finite number of mass points.
\end{abstract}

\begin{keyword} Szeg\H{o} asymptotics \sep orthogonal polynomials \sep matrix-valued measures



\end{keyword}

\end{frontmatter}


\begin{section}{Introduction}

Let $\mu$ be an $l\times l$ matrix-valued Hermitian positive semi-definite finite measure on $\bbR$ of compact support, normalized
by $\mu(\bbR) = \bdone$, where $\bdone$ is the $l\times l$ identity matrix. For any $l\times l$ dimensional matrix functions $f,g$, define
\begin{align}
\lla f,g\rra_{L^2(\mu)}&=\int f(x)^*d\mu(x) g(x);\label{eq1.1}\\
\lla f\rra^2_{L^2(\mu)}&=\lla f,f\rra_{L^2(\mu)}.
\end{align}
Here we can regard $\lla f\rra_{L^2(\mu)}$ as the square root of the non-negative definite matrix $\lla f,f\rra_{L^2(\mu)}$. By $\lla f,g\rra_{L^2}$, with the index just $L^2$, we will mean the product with respect to the Lebesgue measure on the real line or the unit circle, depending on the context.

What we have defined here is the right product of $f$ and $g$, as opposed to the left product $\int f(x)d\mu(x) g(x)^*$, whose properties are completely analogous.

Measure $\mu$ is called non-trivial if $|| \lla f\rra^2_{L^2(\mu)}||>0$ for all polynomials $f$. From now on assume $\mu$ is non-trivial. Then (see \cite{DPS} for the details) there exist unique monic polynomials $P_n$ of degree $n$ satisfying
\begin{equation*}
\lla P_n,f \rra_{L^2(\mu)}=0 \quad \mbox{ for any polynomial } f \mbox{ with } \deg f<n.
\end{equation*}

For any choice of unitary $l\times l$ matrices $\tau_n$ (we demand $\tau_0=\bdone$), the polynomials $p_n=P_n \lla P_n\rra_{L^2(\mu)}^{-1}\tau_n$ are orthonormal:
\begin{equation*}
\lla p_n,p_m \rra_{L^2(\mu)}=\delta_{n,m} \bdone.
\end{equation*}

We are interested in the asymptotic behavior of $p_n$ for the measures $\mu$ whose essential support is a single interval. After scaling and translating, we can assume it is $[-2,2]$:
 \begin{equation}\label{eq1.8}
\esssup\mu=[-2,2].
\end{equation}

Let $\{E_j\}_{j=1}^{N}$ be the point masses of $\mu$ outside $[-2,2]$ counting multiplicities ($N\le \infty$). In~\cite{AN} it was shown that if the absolutely continuous part $f(x)=\frac{d\mu(x)}{dx}$ satisfies the Szeg\H{o} condition
\begin{equation}\label{eq1.9}
\int_{-2}^2 (4-x^2)^{-1/2} \log(\det(f(x)))dx>-\infty,
\end{equation}
and $N$ is finite, then there exists $\lim_{n\to\infty} z^n p_n(z+z^{-1})$ uniformly on the compacts of $\bbD$, and the limit function was constructed more or less explicitly. The scalar case $l=1$ (see~\cite{PY}; another approach is the combination of~\cite{KS} and~\cite{DS1}: see~\cite[Chapter 3]{Rice}) suggests that $N=\infty$ should not really spoil the picture as soon as the condition
\begin{equation}\label{eq1.10}
\sum_{j=1}^N\left(|E_j|-2\right)^{1/2}<\infty
\end{equation}
holds. In fact this condition is necessary if one expects to have the limit $\lim_{n\to\infty} z^n p_n(z+z^{-1})$ to be a Nevanlinna function in $\bbD$.

Assume from now on that~(\ref{eq1.8}),~(\ref{eq1.9}) and~(\ref{eq1.10}) hold. We prove in Theorem~\ref{th2} below that under these assumptions $\lim_{n\to\infty} z^n p_n(z+z^{-1})$ exists uniformly in $\bbD$ and we give a characterization of the limit function. The results are the exact extension to the matrix-valued case of~\cite{PY}, and include \cite[Thm. 2]{AN} as its special case ($N<\infty$).

To prove the result, Aptekarev and Nikishin in~\cite{AN} used an induction on the number of the point masses of $\mu$, which does not work if there are infinitely many of them. The approach used here is similar to the one used in~\cite{PY} for the scalar case (which in turn is an extension of the original Szeg\H{o} proof for the no-bound problem, see~\cite{Sz}). Namely, we first construct a Nevanlinna function $L(z)$ (Section 3), and then consider a certain inner product which, when handled with care, proves that the limit of $z^n p_n(z+z^{-1})$ is indeed $L$ (Section 4).

We start by collecting some auxiliary statements in Section 2.

\bigskip

\textbf{Acknowledgements.} The author would like to thank Barry Simon for helpful discussions.

\end{section}

\begin{section}{Preliminaries}

%\subsection{Matrix OPRL}

\textbf{2.1. Matrix OPRL.}
For the proofs and additional results on the theory of matrix-valued orthogonal polynomials, see the review article of Damanik--Pushnitski--Simon~\cite{DPS}.

Just as in the scalar case, for any choice of unitary $l\times l$ matrices $\tau_n$ ($\tau_0=\bdone$), the orthonormal polynomials $p_n=P_n \lla P_n\rra_{L^2(\mu)}^{-1}\tau_n$ satisfy the recurrence relation
\begin{equation*}
x p_n(x)=p_{n+1}(x)A_{n+1}^*+p_n(x)B_{n+1}+p_{n-1}(x) A_n, \quad n=1,2,\ldots,
\end{equation*}
where $A_n=\lla p_{n-1},xp_n \rra_{L^2(\mu)}$, $B_n=\lla p_{n-1},xp_{n-1} \rra_{L^2(\mu)}$ (with $p_{-1}=\bdnot$, $A_0=\bdone$, the relation holds for $n=0$ too). The corresponding block Jacobi matrix is
\begin{equation*}
J=\left(
\begin{array}{cccc}
B_1&A_1&\mathbf{0}&\cdots\\
A_1^{*}&B_2&A_2&\cdots\\
\mathbf{0}&A_2^*&B_3&\cdots\\
\vdots&\vdots&\vdots&\ddots\end{array}\right).
\end{equation*}

Two block Jacobi matrices $J$ and $\widetilde{J}$ are called equivalent if they correspond to the same measure $\mu$ (but a different choice of $\tau_n$'s). They are equivalent if and only if their Jacobi parameters satisfy
\begin{equation}\label{eq2.3}
\widetilde{A}_n=\sigma_n^* A_n \sigma_{n+1}, \quad
\widetilde{B}_n=\sigma_n^* B_n \sigma_n
\end{equation}
for unitary $\sigma_n$'s with $\sigma_1=\bdone$ (the connection with $\tau_j$'s is $\sigma_n=\tau_{n-1}^* \widetilde{\tau}_{n-1}$). It is easy to see that
\begin{equation}\label{eq2.4}
\widetilde{p}_n(x)=p_n(x)\sigma_{n+1},
\end{equation}
where $\widetilde{p}_n$ are the orthonormal polynomials for $\widetilde{J}$.

We say that a block Jacobi matrix is of type $1$ if $A_n>0$ for all $n$,
of type $2$ if $A_1A_2\ldots A_n>0$ for all $n$, and of type $3$ if every $A_n$ is lower triangular. Each equivalence class of block Jacobi matrices contains exactly one matrix of type $1$, $2$ and $3$.

Define the $m$-function of the measure $\mu$ to be the meromorphic in $\bbC\setminus [-2,2]$ matrix-valued function $m(z)=\int \frac{d\mu(x)}{x-z}$. Define also
\begin{equation*}
M(z)=-m(z+z^{-1}), \quad z\in\bbD.
\end{equation*}
Just as in the scalar case, one easily sees that
\begin{align*}
\imag M(e^{i\theta})&=\pi f(2\cos\theta), \quad 0\le\theta\le\pi,\\
\imag M(e^{i\theta})&=-\pi f(2\cos\theta), \quad -\pi\le\theta\le0,
\end{align*}
where $\imag T\equiv \frac{T-T^*}{2i}$.

Denote
\begin{equation}\label{eq2.5}
\{z_k\}_{k=1}^N=\left\{z\in\bbD \biggm| z=\frac12\left(E_k-\sqrt{{E_k}^2-4}\right)\right\}=\left\{z\in\bbD \mid z+z^{-1}=E_k\right\},
\end{equation}
enumerated in increasing order of their absolute values ($N\le\infty$). Let us assume each $z_k$ is different, and let $m_k$ be the multiplicity of $z_k+z_k^{-1}$ as the eigenvalue. Then~(\ref{eq1.10}) implies
\begin{equation}\label{eq2.6}
\sum_{k=1}^N m_k\log |z_k|<\infty.
\end{equation}

We will be using the so-called $C_0$ Sum Rule from~\cite{DKS}. In a slightly changed form, it looks as follows.
\begin{theorem}[Damanik--Killip--Simon \cite{DKS}]\label{th1}
Suppose $\esssup\mu=[-2,2]$ and $\{z_k\}_{k=1}^N$ be as in~\eqref{eq2.5}. Let
\begin{align*}
Z(J)&= -\frac{1}{2}\int_{-\pi}^\pi \log\det\frac{\imag M(e^{i\theta})}{\sin\theta}\dt,\\
\mathcal{E}_0(J)&= -\sum_{k=1}^N m_k\log |z_k|,\\
A_0(J)&= -\lim_{n\to\infty} \sum_{j=1}^n \log\det |A_j|.
\end{align*}
If any two of $Z, \mathcal{E}_0, A_0$ are finite, then so is the third, and
\begin{equation*}
Z(J)=\mathcal{E}_0(J)+A_0(J).
\end{equation*}
\end{theorem}

\begin{remarks} 1. Here $|T|\equiv \sqrt{T^*T}$.
\smallskip

2. The minus in the expression for $\mathcal{E}_0(J)$ comes from the fact that we chose $z\in\bbD$ in~(\ref{eq2.5}) as opposed to $z\in\bbC\setminus\overline{\bbD}$ in~\cite{DKS}.
\end{remarks}

\smallskip

Note that $z^n p_n(z+z^{-1})$ at $z=0$ is equal to $\kappa_n=\left({{{A}_1}^*}\right)^{-1} \ldots \left({{A}_n}^*\right)^{-1}>0$. Since in our case $\mathcal{E}_0(J)<\infty$ (see~(\ref{eq2.6})) and $Z(J)<\infty$ (follows from~(\ref{eq1.8})), this theorem allows us to find the limit of the determinants of $z^n p_n(z+z^{-1})$ at $z=0$. We will see later that this limit is equal to $\det L(0)$. This will be  helpful in showing that $z^n p_n(z+z^{-1})$ at $z=0$ converges to $L(0)$, and this in turn will allow us to show the uniform convergence in $\bbD$.

%\subsection{Properties of the Matrix Product}
\bigskip

\textbf{2.2. Properties of the Matrix Product.}
We will need a couple of facts about the product we defined in~(\ref{eq1.1}).

\begin{lemma}\label{lm1}  Let $L^2(\bdone\dt)$ be the space of all matrix-valued functions, each entry of which is a scalar $L^2(\dt)$-function.

(a) The following formulae
\begin{align*}
\|f\|_{L^2,1}&\equiv\left(\int_{-\pi}^\pi \|f(\theta)\|^2\dt\right)^{1/2}, \\
\|f\|_{L^2,2}&\equiv\left\|\int_{-\pi}^\pi f(\theta)^* f(\theta)\dt\right\|^{1/2}=|| \lla f,f\rra_{L^2} ||^{1/2}
\end{align*}
define two equivalent (semi)norms on $L^2(\bdone\dt)$:
\begin{equation*}
\|f\|_{L^2,2}\le \|f\|_{L^2,1}\le {l}^{1/2} \|f\|_{L^2,2}.
\end{equation*}

(b) For any $f,g\in L^2$,
\begin{equation*}
\left\|\lla f,g\rra_{L^2}\right\|\le l \|f\|_{L^2,2} \|g\|_{L^2,2}.
\end{equation*}

(c) If $f\in L^2(\bdone\dt)$, then its $n$-th matrix Fourier coefficient $\lla e^{in\theta}I,f\rra_{L^2}\to \mathbf{0}$ as $n\to\infty$.

%(d) If $f_n\to f$, $g_n\to g$ in $L^2(I\dt)$ then
%\begin{equation*}
%\lla f_n,g_n \rra_{L^2} \to \lla f,g \rra_{L^2}
%\end{equation*}

\begin{proof}
(a) The first inequality is obvious. The second follows from
\begin{equation*}
\|f\|_{L^2,1}^2=\int_{-\pi}^\pi \|f(\theta)\|^2\dt\le \int_{-\pi}^\pi \tr(f(\theta)^*f(\theta))\dt =
\tr\left(\int_{-\pi}^\pi f(\theta)^*f(\theta)\dt\right)
\le l \|f\|_{L^2,2}^2.
\end{equation*}

%Now assume $\|f\|^2_2=\left\|\int_{-\pi}^\pi f(\theta)^* f(\theta)\dt\right\|<\infty$. Then in particular for any unit vector $\phi$,
%\begin{equation}
%\int_{-\pi}^\pi \|f(\theta)\phi\|^2 \dt =\left(\int_{-\pi}^\pi f(\theta)^* f(\theta)\dt\,\phi,\phi\right)\le \left\|\int_{-\pi}^\pi f(\theta)^* f(\theta)\dt\right\|
%%\le\tr\int_{-\pi}^\pi f(\theta)^* f(\theta)\dt<\infty
%\end{equation}
%Taking $\phi=e_j$, we obtain $\int_{-\pi}^\pi \sum_{k=1}^l |f_{kj}(\theta)|^2 \dt \le\|f\|_{L^2,2}^2$. Thus
%\begin{equation}
%\|f\|_{L^2,1}^2=\int_{-\pi}^\pi \|f(\theta)\|^2\dt\le \int_{-\pi}^\pi \tr(f(\theta)^*f(\theta))\dt =\int_{-\pi}^\pi \sum_{j,k=1}^l |f_{kj}(\theta)|^2\dt\le l \|f\|_{L^2,2}^2.
%\end{equation}

(b) Using H\"{o}lder, and the equivalence from (a), we get
\begin{equation*}
\left\|\lla f,g\rra_{L^2}\right\|\le
\int_{-\pi}^\pi \|g(\theta)\|\, \|f(\theta)\|\dt \le \|f\|_{L^2,1} \|g\|_{L^2,1}\le l \|f\|_{L^2,2} \|g\|_{L^2,2}.
\end{equation*}

(c) %, (d)%
Follows by looking at each entry separately.
\end{proof}
\end{lemma}

\medskip

%The scalar limit function in~\cite{PY} consisted of an outer factor and a Blaschke product. We will need the matrix-valued analogues of both.

Just as in the scalar case (see~\cite{PY}) we will be expecting the limit function to have the factorized form of a product of a matrix outer factor and a matrix Blaschke product.

%\subsection{Matrix Outer Functions}
\bigskip

\textbf{2.3. Matrix Outer Functions.}
Recall that a scalar analytic function $G$ on $\bbD$ is called outer if it can be recovered from its boundary values $G(e^{i\theta})\equiv\lim_{r\nearrow 1} G(re^{i\theta})$ by the formula
\begin{equation}\label{eq2.18}
G(z)=c \exp\left\{\int_{-\pi}^{\pi} \frac{e^{i\theta}+z}{e^{i\theta}-z} \log |G(e^{i\theta})|\dt\right\}
\end{equation}
for some constant $|c|=1$. Note that it is necessary and sufficient $\log |G(e^{i\theta})|$ to be integrable.

\begin{lemma}[Wiener--Masani \cite{WM}]\label{lm2}
Suppose $w(\theta)$ is a matrix-valued function on the unit circle satisfying
\begin{equation*}
\int_{-\pi}^\pi \log \det w(\theta)\dt>-\infty.
\end{equation*}
Then there exists a unique matrix-valued $H_2(\bbD)$ function $G(z)$ satisfying
\begin{gather}
G(e^{i\theta})^* G(e^{i\theta})=w(\theta),\\
G(0)^*=G(0)>0,\\
\log|\det G(0)|=\int_{-\pi}^{\pi} \log|\det G(e^{i\theta})|\dt.\label{eq2.21}
\end{gather}
\end{lemma}

This is a well-known result of Wiener--Masani \cite{WM}. The proof of the uniqueness part can be found, e.g., in~ \cite{DGK}.

Equality~(\ref{eq2.21}) implies (see \cite[\S 17.17]{Rud}) that $\det G(z)$ is a scalar outer function, which implies (by definition) that $G(z)$ is a matrix-valued outer function. %Though it's irrelevant for the proof,
It follows from \cite[Thm.~2]{Ginz} that there exists a Hermitian matrix-valued integrable function $M(\theta)$ such that
\begin{equation}\label{eq2.23}
\tr \,M(\theta)=\log|\det G(e^{i\theta})|
\end{equation} and
\begin{equation}\label{eq2.24}
G(z)=\rho \stackrel{\curvearrowright}{\int_{-\pi}^{\pi}} \exp\left\{\frac{e^{i\theta}+z}{e^{i\theta}-z} M(\theta)\dt\right\},
\end{equation}
where $\stackrel{\curvearrowright}{\int_{-\pi}^{\pi}}$ is the Potapov multiplicative integral (see~\cite{Pot})
\begin{equation*}
\stackrel{\curvearrowright}{\int_{-\pi}^{\pi}} \exp\left\{F(\theta)\dt\right\}=\lim_{\Delta\theta_j\to0} \stackrel{\curvearrowright}{\prod_{j=0}^{n-1}} e^{F(\phi_j)\Delta\theta_j}, \quad -\pi=\theta_0\le\phi_0\le\theta_1\le\phi_1\le\cdots\le\theta_{n-1}\le\phi_{n-1}\le \theta_n=\pi.
\end{equation*}
The arrow above the product sign simply defines the order of the multiplication in the matrix-valued product. $\rho$ in~(\ref{eq2.24}) is a constant unitary matrix which makes the right-hand side of (\ref{eq2.24}) positive-definite.

Clearly~(\ref{eq2.23})--(\ref{eq2.24}) becomes~(\ref{eq2.18}) if $l=1$.

%\subsection{Blaschke--Potapov Products}
\bigskip

\textbf{2.4. Blaschke--Potapov Products.}
The Blaschke--Potapov elementary factor is a generalization of scalar Blaschke factors (for those familiar with the Potapov theory of $J$--contractive matrix functions: we are considering the signature matrix $J$ to be just the identity matrix $\bdone$):
%$$
%b_n(z,z_n)= U_n^* \left(\left(\frac{|z_j|}{z_j}\frac{z_j-z}{1-z_j z}
%\bdone_{m_j}+(\bdone-\bdone_{m_j})\right)\right)
% U_n,
%$$
$$
B_{z_j,s,U}(z)=U^*\left(
\begin{array}{cccccc}
\frac{|z_j|}{z_j}\frac{z_j-z}{1-z_j z}&0&0&0&\cdots&0\\
0&\ddots&0&0&\cdots&0\\
0&0&\frac{|z_j|}{z_j}\frac{z_j-z}{1-z_j z}&0&\cdots&0\\
0&0&0&1&\cdots&0\\
\vdots&\vdots&\vdots& &\ddots& \\
0&0&0&0&\cdots&1
\end{array}\right)U, \quad z\in\bbD,
$$
where $z_j\in\bbD$, $s$ is the number of the scalar Blaschke factors on the diagonal ($0\le s\le l$), and $U$ is a unitary constant matrix. Clearly $B_{z_j,s,U}$ is an analytic  in $\bbD$ function with unitary values on the unit circle.

The well-known result for the convergence of the scalar Blaschke products is still valid for the matrix-valued case: if
\begin{equation*}
\sum_{k=1}^\infty (1-|z_k|)<\infty,
\end{equation*}
then the product
\begin{equation*}
\stackrel{\curvearrowright}{\prod_{j=1}^\infty} B_{z_j,s_j,U_j}(z)
\end{equation*}
converges uniformly on the compacts of the unit disk (see~\cite{Pot} and~\cite{Ginz2} where this is proven even more generally for the operator-valued setting). The limit function is holomorphic in $\bbD$ with unitary boundary values (see \cite{Arov}).

We have freedom here in the choice of the unitary matrices $U_j$ and numbers $s_j$. We will make use of it in the following lemma.

\begin{lemma}
Let $\{z_k\}_{k=1}^\infty$ with $\sum_{k=1}^\infty (1-|z_k|)<\infty$ be given, with all $z_k$ pairwise different. For any sequence of subspaces $V_k\subset \bbC^l$, there exists a unique product $B(z)=\stackrel{\curvearrowright}{\prod_{j=1}^\infty} B_{z_j,s_j,U_j}(z)$ for some choice of numbers $s_k$, $0\le s_k\le l$, and unitary matrices $U_k$, that  satisfies
\begin{equation}\label{eq2.28}
\ker \res_{z=z_k} B(z)^{-1}=V_k \quad \mbox{for all } k.
\end{equation}
\begin{proof}
Easy induction does the job. Let $I_m$ ($0\le m \le l$) be the diagonal $l\times l$ matrix with first $m$ diagonal elements $1$ and the rest $0$, and $B_n(z)=\stackrel{\curvearrowright}{\prod_{j=1}^n}B_{z_j,s_j,U_j}(z)$ be the partial finite product. Assume that we already chose $\{s_k\}_{k=1}^{n-1}$ and $\{U_k\}_{k=1}^{n-1}$  so that $B_{n-1}(z)$ satisfies
\begin{equation*}
\ker \res_{z=z_k} B_{n-1}(z)^{-1}=V_k, \quad 1\le k \le n-1.
\end{equation*}
Observe that this implies~(\ref{eq2.28}) holds for $1\le k\le n-1$ as well. Put $s_n=l-\dim V_n$. Note that
\begin{equation}\label{eq2.29}
B_{z_n,s_n,U_n}(z)=U_n^* \left(\frac{|z_n|}{z_n}\frac{z_n-z}{1-z z_n}I_{s_n}+(I-I_{s_n})\right) U_n
\end{equation}
and
\begin{equation*}
\ker \res_{z=z_n} B(z)^{-1}=\ker \res_{z=z_n} B_{n}(z)^{-1}=\ker I_{s_n} U_n B_{n-1}(z_n)^{-1}.
\end{equation*}
Note that $B_{n-1}(z_n)$ is invertible (as $z_n\notin\{z_1,\ldots,z_{n-1}\}$), so we can put $U_n$ to be any unitary matrix taking the subspace $B_{n-1}(z_n)^{-1} V_n$ to $\langle \delta_{s_n+1}\cdots\delta_l\rangle$. Note that the choice of $U_n$ is not unique, but the factor $B_{z_n, s_n, U_n}$ is uniquely defined.
\end{proof}
\end{lemma}

\end{section}

\begin{section}{Construction of the limit function}

Let $J$ be the type $2$ Jacobi matrix corresponding to $\mu$, and let $p_n$ be the orthonormal polynomials for $J$.

Let $\nu=\operatorname{Sz}(\mu|_{[-2,2]})$ be the image measure on $\partial\bbD$ of $\mu|_{[-2,2]}$ under  $\theta\mapsto2\cos\theta$: $\int_{-\pi}^\pi g(2\cos\theta)d\nu(\theta)=\int_{-2}^2 g(x) d\mu(x)$ for measurable $g$'s. This is what is called the Szeg\H{o} mapping. Let the Lebesgue decomposition of $\nu$ %with respect to the Lebesgue measure
be
\begin{equation*}
d\nu(\theta)=w(\theta)\frac{d\theta}{2\pi}+d\nu_s.
\end{equation*}
Then
\begin{equation}\label{eq3.2}
w(\theta)=2\pi |\sin\theta| f(2\cos\theta),
\end{equation}
and so~(\ref{eq1.9}) implies
\begin{equation*}\label{eq3.3}
\int_{-\pi}^\pi \log \det w(\theta)\dt>-\infty.
\end{equation*}

Therefore Lemma 2 applies, so there exists a matrix-valued outer $H_2(\bbD)$-function $G(z)$ such that
\begin{gather}
w(\theta)= G (e^{i\theta})^* G(e^{i\theta}), \label{eq3.4}\\
G(0)^*=G(0)>0, \label{eq3.5}\\
\log|\det G(0)|=\int_{-\pi}^{\pi} \log|\det G(e^{i\theta})|\dt. \label{eq3.6}
\end{gather}

Denote $w_k$ to be the weight of $\mu$ at $z_k+z_k^{-1}$:
\begin{equation*}
w_k=\mu(z_k+z_k^{-1}).
\end{equation*}

Now apply Lemma 3 to obtain the Blaschke--Potapov product $B(z)=\stackrel{\curvearrowright}{\prod_{j=1}^\infty} B_{z_j,s_j,U_j}(z)$ (by~(\ref{eq2.6}) it converges) satisfying

\begin{equation}\label{eq24}
\ker \res_{z=z_k} \left(B(z)^{-1}G(z)\right)=\ker w_k \quad \mbox{for all } k.
\end{equation}
Indeed, note that $G(z_k)$ is invertible for any $k$ (since $\det G$ is outer, it can't vanish in $\bbD$), so we can apply Lemma 3 with $V_k=G(z_k) \ker w_k$.

%This can be done inductively. Let $I_m$ ($0\le m \le l$) be the diagonal $l\times l$ matrix with first $m$ diagonal elements $1$ and the rest $0$, and $B_n(z)=\stackrel{\curvearrowright}{\prod_{k=1}^n}b_k(z,z_k)$ be the partial finite product. Assume we defined $B_{n-1}(z)$ satisfying
%\begin{equation}
%\ker \res_{z=z_k} \left(B_{n-1}(z)^{-1}G(z)\right)=\ker w_k, \quad 1\le k \le n-1.
%\end{equation}
%Observe that this implies~(\ref{eq24}) holds for $1\le k\le n-1$ as well. Denote $m_n$ the multiplicity of $z_n+z_n^{-1}$ as an eigenvalue of $J$, i.e. $m_n=l-\dim \ker w_n$. Then put
%\begin{equation}\label{eq2.29}
%b_n(z,z_n)=U_n^* \left(\frac{|z_n|}{z_n}\frac{z-z_n}{1-z z_n}I_{m_n}+(I-I_{m_n})\right) U_n
%\end{equation}
%Then
%\begin{equation}
%\ker \res_{z=z_n} \left(B_{n}(z)^{-1}G(z)\right)=\ker I_{m_n} U_n B_{n-1}(z_n)^{-1}G(z_n)
%\end{equation}
%Observe that $B_{n-1}(z_n)^{-1}G(z_n)$ is invertible ($G$, being outer, doesn't vanish), so we can put $U_n$ to be any unitary matrix taking the subspace $G(z_n)^{-1} B_{n-1}(z_n) \ker w_n$ to $\langle \delta_{m_n+1}\cdots\delta_l\rangle$.

Define for $z\in\bbD$,
%\begin{equation}
%D(z)=U B(z)^{-1}G(z),
%\end{equation}
\begin{equation}\label{eq3.9}
L(z)=\frac{1}{\sqrt2}\,G(z)^{-1}B(z)V,
\end{equation}
where $V$ is a constant unitary such that $L(0)>0$.

Now we can formulate the main result of the paper.

\begin{theorem}\label{th2} Let $\mu$ satisfy~\eqref{eq1.8},~\eqref{eq1.9},~\eqref{eq1.10}. Assume $J$ is of type $2$, and let $\widetilde{J}$ be any equivalent to it matrix with Jacobi parameters~\eqref{eq2.3} and orthonormal polynomials $\widetilde{p}_n$~\eqref{eq2.4}. Assume $\sigma=\lim_{n\to\infty}\sigma_n$ exists. Then
\begin{align}
z^n\widetilde{p}_n\left(z+z^{-1}\right) & \to L(z)\sigma \quad \mbox{uniformly on compacts of } \bbD; \label{eq3.10}\\
\widetilde{p}_n(2\cos\theta) &  =\frac{1}{\sqrt2 } \left(e^{-in\theta}L(e^{i\theta})+e^{in\theta} L(e^{-i\theta}))\right) \sigma +o(1)\quad \mbox{in } L^2\left(w(\theta)\dt\right) \mbox{ sense;} \label{eq3.11}\\
\lla \widetilde{p}_n(x)\rra_{L^2(\mu_s)} & \to \mathbf{0}, \label{eq3.12}
\end{align}
where $w$ is defined in~\eqref{eq3.2}.

The limit function $L$ has a factorization~\eqref{eq3.9}, where $G$ is the unique $H_2(\bbD)$-function satisfying~\eqref{eq3.4}--\eqref{eq3.6} $($and thus has the form~\eqref{eq2.23}--\eqref{eq2.24}$)$, $B$ is a Blaschke--Potapov product, $V$ is a unitary matrix. We have
\begin{align*}
\ker \res_{z=z_k} L(z)^{-1} &= \ker w_k \quad \mbox{for all } k;\\
L(0) &>0.
\end{align*}
\end{theorem}

\medskip

\begin{remarks} 1. We will show that the asymptotics holds for type $2$ Jacobi matrix. Thus by~(\ref{eq2.4}), the polynomials $\widetilde{p}_n$ obey Szeg\H{o} asymptotics if and only if the limit $\lim_{n\to\infty}\sigma_n$ exists, so this condition is also necessary.

\medskip

2. The equivalent way of writing~(\ref{eq3.11}) is
\begin{equation*}
G(e^{i\theta})\widetilde{p}_n(2\cos\theta) = \frac{1}{\sqrt2 } \left(e^{-in\theta}B(e^{i\theta})+e^{in\theta}G(e^{i\theta}) G(e^{-i\theta})^{-1} B(e^{-i\theta})\right)V\sigma +o(1) \quad \mbox{in } L^2\left(\bdone \dt\right) \mbox{ sense.}
\end{equation*}
\end{remarks}

%\medskip
%
%3. The necessity of~(\ref{eq3.13}) was shown in~\cite{AN}
%
\medskip

\noindent Using results from Section $14$ of \cite{DKS} we immediately obtain

\begin{corollary}\label{cor}
Assume the Jacobi parameters of $J$ satisfy
\begin{equation}\label{l2cond}
\sum_{n=1}^\infty \left[ \norm{1-A_n A_n^*}+ \norm{B_n}\right]<\infty.
\end{equation}
Then the associated measure $\mu$ satisfies~\eqref{eq1.8},~\eqref{eq1.9},~\eqref{eq1.10}, and so the conclusions of Theorem~\ref{th2} hold.
\end{corollary}
\begin{remarks}
1. As in Theorem \ref{th2} this establishes Szeg\H{o} asymptotics for the type $2$ Jacobi matrix, and for any equivalent to it matrix for which the limit $\lim_{n\to\infty}\sigma_n$ exists. In \cite{Kozhan2} it's shown that under \eqref{l2cond} this limit does exist matrices of type $1$ and $3$ (or more generally, for any $\widetilde{J}$ the $\widetilde{A}_n$-coefficients of which have eventually only real eigenvalues).
\medskip

2. See also \cite{Kozhan3} for another proof of Corollary \ref{cor}.
\end{remarks}

%\medskip
%
%3. The necessity of~(\ref{eq3.13}) was shown in~\cite{AN}
%

\end{section}
\begin{section}{Proof}
The beginning of the proof follows closely the proof of the Lemma in~\cite{PY}. Denote $$s(e^{i\theta})=G(e^{i\theta}) G(e^{-i\theta})^{-1},$$ and consider the following expression. Expanding the product and using~(\ref{eq3.4}), we get
\begin{multline}\label{eq4.1}
\mathbf{0}\le \llangle G(e^{i\theta})\widetilde{p}_n(2\cos\theta)-\frac{1}{\sqrt2 } \left(e^{-in\theta}B(e^{i\theta})+e^{in\theta}s(e^{i\theta}) B(e^{-i\theta})\right)\rrangle^2_{L^2}+\lla \widetilde{p}_n(x)\rra^2_{L^2(\mu_s)} \\=\int_{-\pi}^\pi \widetilde{p}_n(2\cos\theta)^*w(\theta)\widetilde{p}_n(2\cos\theta)\dt+\lla \widetilde{p}_n(x)\rra^2_{L^2(\mu_s)} +\frac12 \lla e^{-in\theta}B(e^{i\theta})+e^{in\theta}s(e^{i\theta}) B(e^{-i\theta})\rra^2_{L^2} \\-\sqrt2\real \lla G(e^{i\theta})\widetilde{p}_n(2\cos\theta), e^{-in\theta}B(e^{i\theta})+e^{in\theta}s(e^{i\theta}) B(e^{-i\theta})\rra_{L^2},
\end{multline}
where by $\real T$ we mean $\frac{T+T^*}2$.

First of all,
\begin{equation}\label{eq4.2}
\int_{-\pi}^\pi \widetilde{p}_n(2\cos\theta)^*w(\theta)\widetilde{p}_n(2\cos\theta)\dt+\lla \widetilde{p}_n(x)\rra^2_{L^2(\mu_s)} =\lla \widetilde{p}_n(x)\rra^2_{L^2(\mu)}=\bdone.
\end{equation}

Now, observe that \begin{equation*}\label{eq4.3}
\begin{aligned}
s(e^{i\theta})^*s(e^{i\theta})&=G(e^{-i\theta})^{-*} G(e^{i\theta})^* G(e^{i\theta}) G(e^{-i\theta})^{-1}\\ &=G(e^{-i\theta})^{-*} w(\theta) G(e^{-i\theta})^{-1}=G(e^{-i\theta})^{-*} w(-\theta) G(e^{-i\theta})^{-1}=\bdone.
\end{aligned}
\end{equation*}

Thus
\begin{equation}\label{eq4.4}
\begin{aligned}
\frac12 \lla e^{-in\theta}B(e^{i\theta})+e^{in\theta}s(e^{i\theta}) B(e^{-i\theta})\rra^2_{L^2}&=\bdone+\real \lla e^{-in\theta}B(e^{i\theta}),e^{in\theta}s(e^{i\theta}) B(e^{-i\theta})\rra_{L^2} \\&= \bdone+\int_{-\pi}^\pi e^{2in\theta}B(e^{i\theta})^* s(e^{i\theta})B(e^{-i\theta})\dt=\bdone+o(1)
\end{aligned}
\end{equation}
since the function $k(\theta)=B(e^{i\theta})^* s(e^{i\theta})B(e^{-i\theta})$ satisfies $\int_{-\pi}^\pi k(\theta)^* k(\theta)\dt=\bdone$, so by parts (a) and (c) of Lemma~\ref{lm1}, its Fourier coefficients converge to the zero matrix.

Note that for any function $g$ on the unit circle we have
%Analogously to the scalar case, it is straightforward to show that $s(e^{i\theta})G(e^{-i\theta})\widetilde{p}_n(2\cos\theta)=G(e^{i\theta})\widetilde{p}_n(2\cos\theta)$, and then that for any function $g$,
\begin{equation*}\label{eq4.5}
\begin{aligned}
\lla G(e^{i\theta})\widetilde{p}_n(2\cos\theta), s(e^{i\theta})g(e^{-i\theta}) \rra_{L^2}&=
\int_{-\pi}^\pi \widetilde{p}_n(2\cos\theta)^* \,G(e^{i\theta})^* \,G(e^{i\theta}) G(e^{-i\theta})^{-1}g(e^{-i\theta}) \dt\\ &= \int_{-\pi}^\pi \widetilde{p}_n(2\cos\theta)^* \,w(\theta) G(e^{-i\theta})^{-1}g(e^{-i\theta}) \dt\\
&=\int_{-\pi}^\pi \widetilde{p}_n(2\cos\theta)^* \,w(-\theta) G(e^{-i\theta})^{-1}g(e^{-i\theta}) \dt \\
&=\int_{-\pi}^\pi \widetilde{p}_n(2\cos\theta)^* \,G(e^{-i\theta})^* \,g(e^{-i\theta}) \dt \\
&= \int_{-\pi}^\pi \widetilde{p}_n(2\cos\theta)^* \,G(e^{i\theta})^* \,g(e^{i\theta}) \dt=\lla G(e^{i\theta})\widetilde{p}_n(2\cos\theta), g(e^{i\theta}) \rra_{L^2},
\end{aligned}
\end{equation*}
so the third term on the right-hand side of~(\ref{eq4.1}) becomes
\begin{equation}\label{eq4.6}
\real \lla G(e^{i\theta})\widetilde{p}_n(2\cos\theta), e^{-in\theta}B(e^{i\theta})+e^{in\theta}s(e^{i\theta}) B(e^{-i\theta})\rra_{L^2}= 2\real\lla G(e^{i\theta})\widetilde{p}_n(2\cos\theta), e^{-in\theta}B(e^{i\theta})\rra_{L^2}
\end{equation}

\begin{lemma}\label{lm4}
Let $\widetilde{p}_n(x)=r_n x^n +\ldots$ $($in other words, $r_n=({{\widetilde{A}_1}^*})^{-1} \cdots ({\widetilde{A}_n}^*)^{-1}$ $)$. Then $r_n$ are uniformly bounded $($with respect to the operator norm$)$.
\begin{proof}
On the one hand, by~(\ref{eq4.2}),
\begin{equation}\label{eq4.7}
\left\|G(e^{i\theta})\widetilde{p}_n(2\cos\theta)\right\|_{L^2,2}
%\le\lla G(e^{i\theta})\widetilde{p}_n(2\cos\theta)\rra^2_{L^2}
\le \| \lla \widetilde{p}_n(x)\rra_{L^2(\mu)}^2\|^{1/2}= 1.
\end{equation}
On the other, by Lemma~\ref{lm1}(a) and subharmonicity of $\|h(\cdot)\|^2$,
\begin{equation*}\label{eq4.8}
\begin{aligned}
\left\|G(e^{i\theta})\widetilde{p}_n(2\cos\theta)\right\|_{L^2,2}\ge
l^{-1/2} \left\|G(e^{i\theta})\widetilde{p}_n(2\cos\theta)\right\|_{L^2,1}&=
l^{-1/2} \left(\int_{-\pi}^\pi\left\| h(e^{i\theta})\right\|^2\dt\right)^{1/2} \\
&\ge l^{-1/2} \|h(0)\|= l^{-1/2} \|G(0)r_n\|,
\end{aligned}
\end{equation*}
where $h(z)\equiv G(z)\left[ z^n \widetilde{p}_n\left(z+\frac1z\right)\right]$ is analytic in $\bbD$. $G(0)$ is invertible, so $r_n$ are uniformly bounded.
\end{proof}
\end{lemma}

The next lemma will allow us to compute the right-hand side of~(\ref{eq4.6}).

\begin{lemma}\label{lm5} Let $r_n$ be as in the previous lemma. Then
\begin{equation}\label{eq4.9}
\lla e^{-in\theta}B(e^{i\theta}), G(e^{i\theta})\widetilde{p}_n(2\cos\theta)\rra_{L^2}=B(0)^{-1}G(0) r_n+o(1).
\end{equation}
\begin{proof}
The partial products $B_N(z)$ converge to $B(z)$ uniformly on compacts of $\bbD$. This implies that each Fourier coefficient of $B(e^{i\theta})-B_N(e^{i\theta})$ goes to $0$ as $N\to\infty$. Since $\|B\|_{L^2,2}=\|B_N\|_{L^2,2}=1$, weak convergence implies the norm convergence $\| B(e^{i\theta})-B_N(e^{i\theta})\|_{L^2,2}\to 0$. Using~(\ref{eq4.7})  and Lemma~\ref{lm1}(b), we can find $N\in\bbN$ such that
\begin{equation}\label{eq4.10}
\left\|\lla e^{-in\theta}(B(e^{in\theta})-B_N(e^{-in\theta})),G(e^{i\theta})\widetilde{p}_n(2\cos\theta)\rra_{L^2}\right\| \le l \| B(e^{i\theta})-B_N(e^{i\theta})\|_{L^2,2} \|G(e^{i\theta})\widetilde{p}_n(2\cos\theta)\|_{L^2,2}<\epsilon
\end{equation}
holds for any $n\in\bbN$. By Lemma~{\ref{lm4}}, we can also assume that for this $N$,
\begin{equation}\label{eq4.11}
\|B(0)^{-1}G(0) r_n-B_N(0)^{-1}G(0) r_n\|<\epsilon
\end{equation}
also holds for any $n$. Now, $B_N(e^{i\theta})^*=B_N(e^{i\theta})^{-1}$, so
\begin{multline}\label{eq4.12}
\lla e^{-in\theta}B_N(e^{i\theta}), G(e^{i\theta})\widetilde{p}_n(2\cos\theta)\rra_{L^2}= \int_{-\pi}^\pi e^{in\theta}B_N(e^{i\theta})^{-1} G(e^{i\theta})\widetilde{p}_n(2\cos\theta)\dt \\= \int_{\partial\bbD}B_N(z)^{-1} G(z)\widetilde{p}_n\left(z+\frac1z\right)z^n\frac{dz}{2\pi i z}= B_N(0)^{-1}G(0)r_n+\sum_{k=1}^N \res_{z=z_k}\left(B_N(z)^{-1}G(z)\right) \widetilde{p}_n(E_k)z_k^{n-1}
\end{multline}
By the construction,~(\ref{eq24}) holds, which implies $\ker \res_{z=z_k} \left(B_N(z)^{-1}G(z)\right)=\ker w_k=\ker w_k^{1/2}$, which allows us to write $\res_{z=z_k} \left(B_N(z)^{-1}G(z)\right)=S_k w_k^{1/2}$ for some matrix $S_k$. Thus,
\begin{equation}\label{eq4.13}
\left\| \sum_{k=1}^N \res_{z=z_k}\left(B_N(z)^{-1}G(z)\right) \widetilde{p}_n(E_k)z_k^{n-1} \right\| \le \sup_{1\le k \le N} \|S_k\| \sum_{k=1}^N \|w_k^{1/2} \widetilde{p}_n(E_k)\|\, |z_k|^{n-1}.
\end{equation}
But $\|w_k^{1/2} \widetilde{p}_n(E_k)\|=(\|\widetilde{p}_n(E_k)^*w_k \widetilde{p}_n(E_k)\|)^{1/2}\le \|\lla \widetilde{p}_n(x)\rra_{L^2(\mu)}\|=1$. Since $N$ was fixed, this proves that the right-hand side of~(\ref{eq4.13}) goes to $0$ when $n\to\infty$. Combining~(\ref{eq4.10}),~(\ref{eq4.11}),~(\ref{eq4.12}) and~(\ref{eq4.13}), we obtain~~(\ref{eq4.9}).
\end{proof}
\end{lemma}

Now, plugging~(\ref{eq4.2}),~(\ref{eq4.4}),~(\ref{eq4.6}) and~(\ref{eq4.9}) into~(\ref{eq4.1}), we obtain
\begin{equation}\label{eq4.14}
\mathbf{0}\le 2\bdone-2\sqrt2 \real\left(B(0)^{-1}G(0) r_n\right)+o(1).
\end{equation}
Observe that~(\ref{eq4.14}) holds for any initial choice of unitaries $\sigma_n$ in~(\ref{eq2.3}). Let $p_n(x)=\kappa_n x^n +\ldots$ (in other words, $\kappa_n=\left({A_1^*}\right)^{-1} \cdots \left(A_n^*\right)^{-1}>0$). Then~(\ref{eq2.4}) gives $r_n=\kappa_n\sigma_{n+1}$. For each $n$, pick unitary $\sigma_{n+1}$ such that $B(0)^{-1}G(0) r_n=B(0)^{-1}G(0) \kappa_n\sigma_{n+1}$ is positive-definite. Then~(\ref{eq4.14}) gives
\begin{equation}\label{eq4.15}
\sqrt2 \,B(0)^{-1}G(0) \kappa_n\sigma_{n+1}\le \bdone+o(1).
\end{equation}

Denote $H_n\equiv \sqrt2 \,B(0)^{-1}G(0) r_n>\bdnot$. Let $\{\eta^{(n)}_s\}_{s=1}^l$ be the eigenvalues of $H_n$ in non-increasing order. $\eta^{(n)}_s>0$ for any $n,s$. Then~(\ref{eq4.15}) implies \begin{equation}\label{eq4.16}
\limsup_{n\to\infty} \eta^{(n)}_s\le 1
\end{equation}
for each $s=1,\ldots,l$.

On the other hand, let us compute the determinant of $H_n$. By~(\ref{eq2.29}) and~(\ref{eq2.21}),
\begin{equation*}\label{eq4.17}
\log\det{B(0)^{-1}G(0)}=-\sum_{k}m_k\log|z_k|+\int_{-\pi}^{\pi} \log|\det G(e^{i\theta})|\dt,
\end{equation*}
and by Theorem~\ref{th1},
\begin{multline}\label{eq4.18}
\lim_{n\to\infty} \sum_{j=1}^n \log\det |\widetilde{A}_j|+\sum_{k}m_k\log|z_k|= \frac{1}{2}\int_{-\pi}^\pi \log\det\frac{\imag M(e^{i\theta})}{\sin\theta}\dt= \frac{1}{2}\int_{-\pi}^\pi \log\det\frac{\pi f(2\cos\theta)}{|\sin\theta|}\dt \\
\begin{aligned}=\frac{1}{2}\int_{-\pi}^\pi \log\det\frac{w(\theta)}{2\sin^2\theta}\dt&=\int_{-\pi}^{\pi} \log|\det G(e^{i\theta})|\dt-\frac{l}{4\pi}\int_{-\pi}^\pi \log\left(2\sin^2\theta\right)d\theta \\&=\int_{-\pi}^{\pi} \log|\det G(e^{i\theta})|\dt+\frac{l}{2} \log2.
\end{aligned}
\end{multline}

Now note that
\begin{equation*}\label{eq4.19}
\log\det r_n= -\sum_{j=1}^n \log\det \widetilde{A}_j^* = -\sum_{j=1}^n \log\det |\widetilde{A}_j|+\log\det \rho_n
\end{equation*}
for some unitary matrix $\rho_n$. However, $H_n>\bdnot$, so $\det \rho_n$ must be $1$ as otherwise $\log\det H_n=\frac{l}{2}\log2+\log\det{B(0)^{-1}G(0)}+\log\det r_n$ cannot be real. Thus we obtain
\begin{multline*}\label{eq4.20}
\begin{aligned}
\log\det H_n&=\frac{l}{2}\log2+\log\det{B(0)^{-1}G(0)}+\log\det r_n\\ &=\frac{l}{2}\log2-\sum_{k}m_k\log|z_k|+\int_{-\pi}^{\pi} \log|\det G(e^{i\theta})|\dt -\sum_{j=1}^n \log\det |\widetilde{A}_j| \to 0
\end{aligned}
\end{multline*}
by~(\ref{eq4.18}). Thus $\lim_{n\to\infty} \det H_n=1$. Together with~(\ref{eq4.16}) this implies $\lim_{n\to\infty} \eta^{(n)}_s= 1$ for each $s$, and so $H_n\to \bdone$. This proves $\kappa_n\sigma_{n+1}\to 2^{-1/2}G(0)^{-1}B(0)$. But $|\kappa_n\sigma_{n+1}|=\kappa_n$ (here temporarily $|T|\equiv \sqrt{TT^*}$ instead of $\sqrt{T^*T}$), so
\begin{equation*}
\kappa_n\to 2^{-1/2}\left|G(0)^{-1}B(0)\right|=L(0).
\end{equation*}
Also, $\sigma_n\to V^*$.

Thus for the chosen $\sigma$'s, the right-hand side of~(\ref{eq4.14}) goes to the zero matrix. This implies
\begin{equation*}
\left\| G(e^{i\theta})\widetilde{p}_n(2\cos\theta)-\frac{1}{\sqrt2 } \left(e^{-in\theta}B(e^{i\theta})+e^{in\theta}s(e^{i\theta}) B(e^{-i\theta})\right)\right\|_{L^2,2}\to 0
\end{equation*}
and
\begin{equation*}
\lla \widetilde{p}_n(x)\rra^2_{L^2(\mu_s)}\to\mathbf{0}.
\end{equation*}
Taking into account that $p_n(x)=\widetilde{p}_n(x)\sigma_{n+1}^*$ and $\sigma_n\to V^*$, we get
\begin{equation*}
\left\| G(e^{i\theta})p_n(2\cos\theta)-\frac{1}{\sqrt2 } \left(e^{-in\theta}B(e^{i\theta})+e^{in\theta}s(e^{i\theta}) B(e^{-i\theta})\right)V\right\|_{L^2,2}\to 0
\end{equation*}
and
\begin{equation*}
\lla p_n(x)\rra^2_{L^2(\mu_s)}\to\mathbf{0}.
\end{equation*}

This proves (\ref{eq3.11})--(\ref{eq3.12}) for the type $2$ case. To prove~(\ref{eq3.10}), by Lemma~\ref{lm1}(b),
\begin{multline}\label{eq4.26}
\left\|\llangle \frac{e^{-in\theta}}{1-e^{i\theta}\bar{z}}\bdone,G(e^{i\theta})p_n(2\cos\theta)-\frac{1}{\sqrt2 } \left(e^{-in\theta}B(e^{i\theta})+e^{in\theta}s(e^{i\theta}) B(e^{-i\theta})\right)V\rrangle_{L^2}\right\| \\
\le l %\left\|\frac{e^{-in\theta}}{1-e^{i\theta}\bar{z}}I\right\|_2
\frac{1}{\sqrt{1-|z|^2}}
\left\|G(e^{i\theta})p_n(2\cos\theta)-\frac{1}{\sqrt2 } \left(e^{-in\theta}B(e^{i\theta})+e^{in\theta}s(e^{i\theta}) B(e^{-i\theta})\right)V\right\|_{L^2,2}\to 0
\end{multline}
uniformly on compacts of $\bbD$.
On the other hand,
\begin{multline}\label{eq4.27}
\llangle \frac{e^{-in\theta}}{1-e^{i\theta}\bar{z}}\bdone,G(e^{i\theta})p_n(2\cos\theta)-\frac{1}{\sqrt2 } e^{-in\theta}B(e^{i\theta})V\rrangle_{L^2} \\= \int_{-\pi}^\pi \frac{e^{in\theta}}{1-e^{-i\theta}z}\left(G(e^{i\theta})p_n(2\cos\theta)-\frac{1}{\sqrt2 } e^{-in\theta}B(e^{i\theta})V\right) \dt
%=\int_{\partial\D} \frac{1}{1-z\zeta^{-1}}
%\left(\zeta^n G(\zeta)P_n(\zeta+\zeta^{-1})
%-\frac{1}{\sqrt2} B(\zeta)U^*\right)
%\frac{d\zeta}{2\pi i \zeta}
=z^n G(z)p_n(z+z^{-1})
-\frac{1}{\sqrt2} B(z)V,
\end{multline}
and
\begin{equation}\label{eq4.28}
\llangle \frac{e^{-in\theta}}{1-e^{i\theta}\bar{z}}\bdone,\frac{1}{\sqrt2 } \left(e^{in\theta}s(e^{i\theta}) B(e^{-i\theta})\right)V\rrangle_{L^2}\to \mathbf{0}\quad \mbox{uniformly on compacts of } \bbD
\end{equation}
by Lemma~\ref{lm1}(c). %(uniformity follows from matrix-valued Parseval's equality).
Together, (\ref{eq4.26}), (\ref{eq4.27}) and (\ref{eq4.28}) give
\begin{equation*}
z^n p_n\left(z+z^{-1}\right)  \to L(z) \quad \mbox{uniformly on compacts of } \bbD.
\end{equation*}
Thus we proved~(\ref{eq3.10})--(\ref{eq3.12}) for the type $2$ case. The result for any equivalent $\widetilde{J}$ with Jacobi parameters~(\ref{eq2.3}) for which the limit $\sigma=\lim_{n\to\infty}\sigma_n$ exists, follows immediately from $\widetilde{p}_n(x)=p_n(x)\sigma_{n+1}$.

\end{section}

\end{document}